\documentclass[11pt,reqno]{amsart}
\usepackage{amsmath,amsfonts, amssymb, amsthm}
  \usepackage{paralist}
  \usepackage{graphics} 
  \usepackage{epsfig} 
\usepackage{graphicx}  \usepackage{epstopdf}
 \usepackage[colorlinks=true]{hyperref}
\hypersetup{urlcolor=blue, citecolor=red}

\newtheorem{theorem}{Theorem}[section]

\newtheorem{lemma}[theorem]{Lemma}

\newtheorem{conjecture}{Conjecture}[section]

\theoremstyle{definition}

\newtheorem{remark}{Remark}[section]

\def\pmod #1{\ ({\rm{mod}}\ #1)}
\def\Z{\Bbb Z}
\def\N{\Bbb N}

\def\Q{\Bbb Q}

\def\l{\left}
\def\r{\right}
\def\bg{\bigg}
\def\({\bg(}
\def\){\bg)}
\def\t{\text}
\def\f{\frac}
\def\mo{{\rm{mod}\ }}
\def\pmod#1{\ (\mo\ #1)}

\def\ls{\leq}
\def\gs{\geq}

\def\bi{\binom}

\def\eq{\equiv}

\def\Proof{\noindent{\it Proof}}

\begin{document}
\hbox{Accepted by Colloq. Math. See also {\tt arXiv:2009.04379}}
\medskip

\title[Some new series for $1/\pi$ motivated by congruences]
      {Some new series for $1/\pi$\\ motivated by congruences}
\author[Zhi-Wei Sun]{Zhi-Wei Sun}


\address{Department of Mathematics, Nanjing
University, Nanjing 210093, People's Republic of China}
\email{zwsun@nju.edu.cn}

\keywords{Ramanujan-type series for $1/\pi$, congruences, binomial coefficients, symbolic computation.
\newline \indent 2020 {\it Mathematics Subject Classification}. Primary 11B65, 05A19; Secondary 11A07, 11E25, 33F10.
\newline \indent Supported by the Natural Science Foundation of China (grant no. 11971222).}

\begin{abstract}
In this paper, we deduce a family of six new series for $1/\pi$; for example,
$$\sum_{n=0}^\infty\frac{41673840n+4777111}{5780^n}W_n\left(\frac{1444}{1445}\right)
=\frac{147758475}{\sqrt{95}\,\pi}$$
where $W_n(x)=\sum_{k=0}^n\binom nk\binom{n+k}k\binom{2k}k\binom{2(n-k)}{n-k}x^k$.
To do so, we manage to transform our series to series of the type
$$\sum_{n=0}^\infty\frac{an+b}{m^n}\sum_{k=0}^n\binom nk^4$$
studied by Shaun Cooper in 2012.
In addition, we pose $17$ new series for $1/\pi$ motivated by congruences; for example, we conjecture that
$$\sum_{k=0}^\infty\frac{4290k+367}{3136^k}\binom{2k}kT_k(14,1)T_k(17,16)=\frac{5390}{\pi},$$
where $T_k(b,c)$ is the coefficient of $x^k$ in the expansion of $(x^2+bx+c)^k$.
\end{abstract}
\maketitle

\section{Introduction}

Let $n\in\N=\{0,1,2,\ldots\}$. In 1894 J. Franel \cite{Fr} introduced the usual Franel numbers $f_n=\sum_{k=0}^n\bi nk^3\ (n\in\N)$ and the Franel numbers
$f_n^{(4)}=\sum_{k=0}^n\bi nk^4\ (n\in\N)$ of order four. By Zeilberger's algorithm (cf. \cite{PWZ}),
 the sequence $(f_n^{(4)})_{n\gs0}$ satisfies the following recurrence first claimed by Franel:
 $$(n+2)^3f_{n+2}^{(4)}=4(1+n)(3+4n)(5+4n)f_n^{(4)}+2(3+2n)(7+9n+3n^2)f_{n+1}^{(4)}.$$

M. Rogers and A. Straub \cite{RS} confirmed the author's conjectural series for $1/\pi$ involving Franel polynomials.

  In 2005 Y. Yang used modular forms of level $10$ to discover the following curious identity relating Franel numbers of order four to Ramanujan-type series for $1/\pi$:
 $$\sum_{k=0}^\infty\f{4k+1}{36^k}f_k^{(4)}=\f{18}{\sqrt{15}\,\pi}.$$
 This has not been published by Yang, but more identities of this kind were deduced by S. Cooper \cite{Co} in 2012 via modular forms.
 For the classical Ramanujan-type series for $1/\pi$, one may consult \cite{BB,Be,R} and the nice survey
 given by Cooper \cite[Chapter 14]{Co17}.

For $n\in\N$ the polynomial
 \begin{align*}W_n(x)=&\sum_{k=0}^n\bi nk\bi{n+k}k\bi{2k}k\bi{2(n-k)}{n-k}x^k
 \\=&\sum_{k=0}^n\bi{n+k}{2k}\bi{2k}k^2\bi{2(n-k)}{n-k}x^k
 \end{align*}
 at $x=-1$ coincides with $(-1)^nf_n^{(4)}$, this can be easily verified since the sequence $((-1)^nW_n(-1))_{n\gs0}$
 satisfies the same recurrence as $(f_n^{(4)})_{n\gs0}$.
In 2011 the author \cite[(3.1)-(3.10)]{S-11} proposed ten identities of the form
 $$\sum_{k=0}^\infty\f{ak+b}{m^k}W_k\l(\f1m\r)=\f{C}{\pi},$$
 where $a,b,m$ are integers with $am\not=0$, and $C^2$ is rational.
 They were later confirmed in \cite{CWZ}.

 In this paper we establish six new series for $1/\pi$ involving $W_n(x)$.

 \begin{theorem}\label{Th1.1} We have the following identities:
 \begin{align}\label{W2}\sum_{k=0}^\infty\f{45k+8}{40^k}W_k\l(\f9{10}\r)&=\f{215\sqrt{15}}{12\pi},
 \\\label{W3}\sum_{k=0}^\infty\f{1360k+389}{(-60)^k}W_k\l(\f{16}{15}\r)&=\f{205\sqrt{15}}{\pi},
\\\label{W6}\sum_{k=0}^\infty\f{735k+124}{200^k}W_k\l(\f{49}{50}\r)&=\f{10125\sqrt7}{56\pi},
\\\label{W8}\sum_{k=0}^\infty\f{376380k+69727}{(-320)^k}W_k\l(\f{81}{80}\r)&=\f{260480\sqrt5}{3\pi},
 \\\label{W12}\sum_{k=0}^\infty\f{348840k+47461}{1300^k}W_k\l(\f{324}{325}\r)&=\f{1314625\sqrt2}{12\pi},
 \\\label{W15}\sum_{k=0}^\infty\f{41673840k+4777111}{5780^k}W_k\l(\f{1444}{1445}\r)
 &=\f{147758475}{\sqrt{95}\,\pi}.
 \end{align}
 \end{theorem}

We also have $9$ conjectural series for $1/\pi$ involving $W_n(x)$ as listed in the following conjecture.

 \begin{conjecture}\label{Conj-W} We have the following identities:
 \begin{align}
 \label{W1}\sum_{k=0}^\infty\f{4k+1}{6^k}W_k\l(-\f18\r)&=\f{\sqrt{72+42\sqrt3}}{\pi},
 \\\label{W4}\sum_{k=0}^\infty\f{392k+65}{(-108)^k}W_k\l(-\f{49}{12}\r)&=\f{387\sqrt3}{\pi},
 \\\label{W5}\sum_{k=0}^\infty\f{168k+23}{112^k}W_k\l(\f{63}{16}\r)&=\f{1652\sqrt3}{9\pi},
  \\\label{W7}\sum_{k=0}^\infty\f{1512k+257}{(-320)^k}W_k\l(-\f{405}{64}\r)&=\f{1184\sqrt{35}}{5\pi},
 \\\label{W9}\sum_{k=0}^\infty\f{56k+9}{324^k}W_k\l(\f{25}4\r)&=\f{1134\sqrt{35}}{125\pi},
 \\\label{W10}\sum_{k=0}^\infty\f{13000k-1811}{(-1296)^k}W_k\l(-\f{625}9\r)&=\f{49356\sqrt{39}}{5\pi},
 \\\label{W11}\sum_{k=0}^\infty\f{9360k-1343}{1300^k}W_k\l(\f{900}{13}\r)&=\f{21515\sqrt{39}}{3\pi},
 \\\label{W13}\sum_{k=0}^\infty\f{56355k+2443}{(-5776)^k}W_k\l(-\f{83521}{361}\r)&=\f{4669535\sqrt2}{68\pi},
 \\\label{W14}\sum_{k=0}^\infty\f{5928k+253}{5780^k}W_k\l(\f{1156}5\r)&=\f{28951\sqrt2}{4\pi}.
  \end{align}
 \end{conjecture}
\begin{remark}\label{R-W} Motivated by congruences, the author actually found \eqref{W2}-\eqref{W14} in 2020.
\end{remark}

Van Hamme \cite{vH} thought that classical Ramanujan-type series for $1/\pi$
should have their $p$-adic analogues involving the $p$-adic Gamma function. This does not hold in general for generalized Ramanujan-type series, for example, the author \cite[Conjecture 1.5]{S13d} discovered the identity
$$\sum_{n=0}^\infty\f{6n-1}{256^n}\bi{2n}n\sum_{k=0}^n\bi{2k}k^2\bi{2(n-k)}{n-k}12^{n-k}=\f{8\sqrt3}{\pi}$$
(which was later confirmed in \cite{CWZ}) and conjectured its related $p$-adic congruence
$$\sum_{n=0}^{p-1}\f{6n-1}{256^n}\bi{2n}n\sum_{k=0}^n\bi{2k}k^2\bi{2(n-k)}{n-k}12^{n-k}\eq-p\pmod{p^2}$$
(with $p$ any prime greater than $3$) which has nothing to do with the Legendre symbol $(\f{-3}p)$.

For the author's philosophy to generate series for $1/\pi$ via congruences, one may consult
the survey \cite{S13d} and the recent paper \cite[Section 1]{S20}.

The so-called ``holonomic alchemy" (cf. \cite{CWZ}) does not work for proving our Theorem \ref{Th1.1},
for the reason see Lemma \ref{Lem-2.1} and Remark \ref{Rem-2.1}.
We will prove Theorem \ref{Th1.1} in the next section via transforming \eqref{W2}-\eqref{W15} to series of the type $$\sum_{k=0}^\infty\f{ak+b}{m^k}f_k^{(4)}$$ studied by Cooper \cite{Co},
and present related conjectural congruences in Section 3.

In Sections 4 and 5, we will pose $10$ other new conjectural series for $1/\pi$ motivated by congruences.

\section{Proof of Theorem \ref{Th1.1}}
 \setcounter{equation}{0}
 \setcounter{conjecture}{0}
 \setcounter{theorem}{0}

 \begin{lemma} \label{Lem-2.1} For $|z|\ls1/30$, we have
 \begin{equation}\label{W_k}\sum_{k=0}^\infty\f{z^k}{(1+4z)^{k+1}}W_k\l(\f1{1+4z}\r)=\sum_{n=0}^\infty f_n^{(4)}z^n
 \end{equation}
 and
 \begin{equation}\label{kW_k}\sum_{k=0}^\infty\f{kz^k}{(1+4z)^{k+1}}W_k\l(\f1{1+4z}\r)=\sum_{n=0}^\infty n(f_n^{(4)}+4s_n)z^n,
 \end{equation}
 where
 \begin{equation}\label{sn}s_n:=\sum_{0\ls j<n}(-1)^{n-1-j}\bi{n-1}j\bi{n+j}j\bi{2j}j\bi{2(n-1-j)}{n-1-j}.
 \end{equation}
 \end{lemma}
 \begin{remark}\label{Rem-2.1} Note that the identity \eqref{kW_k}
 contains a sophisticated term $s_n$ defined by \eqref{sn}.
 It is difficult to see how $s_n$ is related to the Franel numbers of order $4$.
 This is why the ``holonomic alchemy" (cf. \cite{CWZ}) is not helpful to our proof of Theorem \ref{Th1.1}.
 \end{remark}
 \medskip

 \noindent{\it Proof of Lemma \ref{Lem-2.1}}.
  Let $N$ be any nonnegative integer. Then
 \begin{align*}&\sum_{k=0}^N\f{z^k}{(1+4z)^{k+1}}W_k\l(\f1{4z+1}\r)
 \\=&\sum_{k=0}^Nz^k\sum_{j=0}^k\bi{k+j}{2j}\bi{2j}j^2\bi{2(k-j)}{k-j}(1+4z)^{-j-k-1}
 \\=&\sum_{k=0}^Nz^k\sum_{j=0}^k\bi{k+j}{2j}\bi{2j}j^2\bi{2(k-j)}{k-j}\sum_{r=0}^\infty\bi{-j-k-1}r(4z)^r
 \end{align*}
 and hence
 \begin{align*}&\sum_{k=0}^N\f{z^k}{(1+4z)^{k+1}}W_k\l(\f1{4z+1}\r)
 \\=&\sum_{n=0}^\infty z^n\sum_{k=0}^{\min\{n,N\}}\sum_{j=0}^k\bi{k+j}{2j}\bi{2j}j^2\bi{2(k-j)}{k-j}\bi{-j-k-1}{n-k}4^{n-k}
 \\=&\sum_{n=0}^\infty z^n\sum_{j=0}^{\min\{n,N\}}\bi{2j}j^2\sum_{k=j}^{\min\{n,N\}}\bi{k+j}{2j}\bi{2(k-j)}{k-j}\bi{n+j}{k+j}(-4)^{n-k}
\\=&\sum_{n=0}^\infty z^n\sum_{j=0}^{\min\{n,N\}}(-4)^{n-j}\bi{2j}j^2\bi{n+j}{2j}\sum_{k=j}^{\min\{n,N\}}\bi{n-j}{k-j}\f{\bi{2(k-j)}{k-j}}{(-4)^{k-j}}.
\end{align*}
Similarly,
\begin{align*}&\sum_{k=0}^N\f{kz^k}{(1+4z)^{k+1}}W_k\l(\f1{4z+1}\r)
\\=&\sum_{n=0}^\infty z^n\sum_{j=0}^{\min\{n,N\}}(-4)^{n-j}\bi{2j}j^2\bi{n+j}{2j}\sum_{k=j}^{\min\{n,N\}}k\bi{n-j}{n-k}\f{\bi{2(k-j)}{k-j}}{(-4)^{k-j}}.
\end{align*}
Clearly $\bi{2m}m\ls(1+1)^{2m}=4^m$ for all $m\in\N$. Thus
$$\bg|\sum_{k=j}^{\min\{n,N\}}\bi{n-j}{k-j}\f{\bi{2(k-j)}{k-j}}{(-4)^{k-j}}\bg|\ls \sum_{k\gs j}\bi{n-j}{k-j}=2^{n-j}$$
and hence
\begin{align*}&\bg|\sum_{j=0}^{\min\{n,N\}}(-4)^{n-j}\bi{2j}j^2\bi{n+j}{2j}
\sum_{k=j}^{\min\{n,N\}}\bi{n-j}{k-j}\f{\bi{2(k-j)}{k-j}}{(-4)^{k-j}}\bg|
\\\ls&\sum_{j=0}^{\min\{n,N\}}4^n\bi{n+j}{2j}\bi{2j}j2^{n-j}
\ls 8^n\sum_{j=0}^n\bi nj\bi{n+j}j\l(\f{2-1}2\r)^j=8^nP_n(2),
\end{align*}
where
$$P_n(x)=\sum_{k=0}^n\bi nk\bi{n+k}k\l(\f{x-1}2\r)^k$$
is the Legendre polynomial of degree $n$.
Similarly,
\begin{align*}&\bg|\sum_{j=0}^{\min\{n,N\}}(-4)^{n-j}\bi{2j}j^2\bi{n+j}{2j}
\sum_{k=j}^{\min\{n,N\}}k\bi{n-j}{n-k}\f{\bi{2(k-j)}{k-j}}{(-4)^{k-j}}\bg|
\\\ls&\sum_{j=0}^{\min\{n,N\}}4^n\bi{n+j}{2j}\bi{2j}j\min\{n,N\}2^{n-j}
\ls n8^nP_n(2).
\end{align*}
By the Laplace-Heine formula (cf. \cite[p.\,194]{Sz}),
$$P_n(2)\sim \f{(2+\sqrt3)^{n+1/2}}{\sqrt{2n\pi}\root 4\of 3}\ \ \t{as}\ n\to+\infty.$$
As $8(2+\sqrt3)<29.86$, we have $n8^nP_n(2)<30^n$ if $n$ is sufficiently large. Recall that $|z|<1/30$.

In view of the above,
\begin{align*}&\lim_{N\to+\infty}\sum_{k=0}^N\f{z^k}{(1+4z)^k}W_k\l(\f1{1+4z}\r)
\\=&\lim_{N\to+\infty}\sum_{n=0}^Nz^n\sum_{j=0}^n (-4)^{n-j}\bi{2j}j^2\bi{n+j}{2j}
\sum_{k=j}^{n}\bi{n-j}{n-k}\bi{-1/2}{k-j}
\\=&\sum_{n=0}^\infty z^n\sum_{j=0}^n(-4)^{n-j}\bi{2j}j^2\bi{n+j}{2j}\bi{n-j-1/2}{n-j}
\\=&\sum_{n=0}^\infty z^n\sum_{j=0}^n\bi{n+j}{2j}\bi{2j}j^2\bi{2(n-j)}{n-j}(-1)^{n-j}
\\=&\sum_{n=0}^\infty z^n(-1)^nW_n(-1)=\sum_{n=0}^\infty f_n^{(4)}z^n.
\end{align*}
Similarly,
\begin{align*}&\lim_{N\to+\infty}\sum_{k=0}^N\f{kz^k}{(1+4z)^k}W_k\l(\f1{1+4z}\r)-\sum_{n=0}^\infty nf_n^{(4)}z^n
\\=&\lim_{N\to+\infty}\sum_{n=0}^Nz^n\sum_{j=0}^n (-4)^{n-j}\bi{2j}j^2\bi{n+j}{2j}
\sum_{k=j}^{n}(k-n)\bi{n-j}{n-k}\bi{-1/2}{k-j}
\\=&-\sum_{n=0}^\infty z^n
\sum_{j=0}^n (-4)^{n-j}\bi{2j}j^2\bi{n+j}{2j}(n-j)\sum_{j\ls k<n}\bi{n-j-1}{n-k-1}\bi{-1/2}{k-j}
\\=&-\sum_{n=0}^\infty z^n
\sum_{j=0}^n (-4)^{n-j}(n-j)\bi nj\bi{n+j}j\bi{2j}j\bi{n-j-3/2}{n-j-1}
\\=&\sum_{n=0}^\infty nz^n
\sum_{0\ls j<n} 4^{n-j}\bi {n-1}j\bi{n+j}j\bi{2j}j\bi{-1/2}{n-j-1}
\\=&\sum_{n=0}^\infty nz^n
\sum_{0\ls j<n} (-1)^{n-j-1}4\bi {n-1}j\bi{n+j}j\bi{2j}j\bi{2(n-j-1)}{n-j-1}.
\end{align*}
So we have the desired result. \qed

\begin{lemma}\label{Lem2.2} For any $n\in\N$ we have
\begin{equation}\label{s-rec}\begin{aligned}&5n(4n+1)((n+2)s_{n+2}-16ns_n)
\\=&\ (30n^3+54n^2+7n-2)f_{n+1}^{(4)}+2(60n^3+58n^2+17n+2)f_n^{(4)}.
\end{aligned}\end{equation}
\end{lemma}
\Proof. Let $u_n$ denote the left-hand side or the right-hand side of \eqref{s-rec}. Via Zeilberger's algorithm, we find that
\begin{align*}&(1+n) (3+n)^3 (5+4 n)u_{n+2}
\\&\times (344+2572 n+8198 n^2+13329 n^3+10875 n^4+4190 n^5+600 n^6)
\\=&\ 2 (2+n) (9+4 n)P(n)u_{n+1}+4(1+n) (2+n) (3+4 n) (5+4 n) (9+4 n)Q(n) u_n
\end{align*}
for all $n=0,1,2,\ldots$, where
\begin{align*}P(n)=&\ 62208+506208 n+1799416 n^2+3578972 n^3+4250502 n^4
\\&\ +3104119 n^5+1401609 n^6+380700 n^7+56940 n^8+3600 n^9
\end{align*}
and
$$Q(n)=40108+127005 n+164335 n^2+110729 n^3+40825 n^4+7790 n^5+600 n^6.$$
Note also that $u_0=0$, $u_1=2150$ and $u_2=103680$. As both sides of \eqref{s-rec}
give the same integer sequence $(u_n)_{n\gs0}$, we have \eqref{s-rec} as desired. \qed

Now we are able to present an auxiliary theorem.

\begin{theorem}\label{auxi} Let $a,b$ and $x$ be complex numbers with $|x-1|\gs7.5$. Then
\begin{equation}\label{x-iden}\begin{aligned}&\f {10}x(x-1)^2(x-2)\sum_{n=0}^\infty\f{an+b}{(4x)^n}W_n\l(1-\f1x\r)
\\=&\sum_{k=0}^\infty
(2ax(5x-7)k+a(10x-13)+10b(x-1)(x-2))\f{f_k^{(4)}}{(4x-4)^k}.
\end{aligned}
\end{equation}
\end{theorem}
\Proof. Note that $|1/(4x-4)|\ls1/30$. Applying \eqref{W_k} with $z=1/(4x-4)$, we get
\begin{equation}\label{W_k-x}\sum_{n=0}^\infty\f1{(4x)^n}W_n\l(1-\f1x\r)=\f x{x-1}\sum_{k=0}^\infty\f{f_k^{(4)}}{(4x-4)^k}.
\end{equation}
If we have
\begin{equation}\label{kW_k-x}
\sum_{n=0}^\infty\f{n}{(4x)^n}W_n\l(1-\f1x\r)=\f x{10(x-1)^2(x-2)}\sum_{k=0}^\infty\f{(10x-14)(kx+1)+1}{(4x-4)^k}f_k^{(4)},
\end{equation}
then combining \eqref{W_k-x} with \eqref{kW_k-x} we immediately get \eqref{x-iden}.
The identity \eqref{kW_k-x} is equivalent to the following one with $z=1/(4x-4)$:
\begin{equation}\label{z-iden}\begin{aligned}&5(1-4z)\sum_{k=0}^\infty\f{kz^k}{(1+4z)^{k+1}}W_k\l(\f1{1+4z}\r)
\\=&\sum_{k=0}^\infty((5-8z)(1+4z)k+4z(5-6z))f_k^{(4)}z^k.
\end{aligned}
\end{equation}

Below we prove \eqref{z-iden} for $|z|\ls 1/30$.
For convenience, we write $[z^m]f(z)$ with $m\in\N$ to denote the coefficient of $z^m$ in the power series expansion of $f(z)$.

By Lemma 2.1, for any $n\in\N$ we have
\begin{align*}&[z^{n+1}](1-16z^2)\sum_{k=1}^\infty\f{kz^{k-1}}{(1+4z)^{k+1}}W_k\l(\f1{1+4z}\r)
\\=\ &[z^{n+2}](1-16z^2)\sum_{m=0}^\infty m(f_m^{(4)}+4s_m)z^m
\\=\ &(n+2)(f_{n+2}^{(4)}+4s_{n+2})-16n(f_n^{(4)}+4s_n)
\\=\ &(n+2)f_{n+2}^{(4)}-16nf_n^{(4)}+4((n+2)s_{n+2}-16ns_n).
\end{align*}

Now let $n\in\Z^+$. By the recurrence of $(f_m^{(4)})_{m\gs0}$, we have
$$4n(4n+1)(4n-1)f_{n-1}^{(4)}=(n+1)^3f_{n+1}^{(4)}-2(2n+1)(3n^2+3n+1)f_n^{(4)}$$
and hence
\begin{align*}&n(4n+1)((32n+52)f_{n+1}^{(4)}+(96n+56)f_n^{(4)}-32(4n-1)f_{n-1}^{(4)})
\\=&\ 4n(4n+1)(8n+13)f_{n+1}^{(4)}+8n(4n+1)(12n+7)f_n^{(4)}
\\&\ -8(n+1)^3f_{n+1}^{(4)}+16(2n+1)(3n^2+3n+1)f_n^{(4)}
\\=&\ 4(30n^3+54n^2+7n-2)f_{n+1}^{(4)}+8(60n^3+58n^2+17n+2)f_n^{(4)}
\\=&\ 20n(4n+1)((n+2)s_{n+2}-16ns_n)
\end{align*}
with the aid of Lemma \ref{Lem2.2}.
Combining this with the last paragraph, we get
\begin{align*}&[z^{n+1}]5(16z^2-1)\sum_{k=1}^\infty\f{kz^{k-1}}{(1+4z)^{k+1}}W_k\l(\f1{1+4z}\r)
\\=\ &-5(n+2)f_{n+2}^{(4)}+80nf_{n}^{(4)}-20((n+2)s_{n+2}-16ns_n)
\\=\ &-5(n+2)f_{n+2}^{(4)}+80nf_{n}^{(4)}
-(32n+52)f_{n+1}^{(4)}
\\&\ -(96n+56)f_n^{(4)}+32(4n-1)f_{n-1}^{(4)}
\\=&\ [z^{n+1}](32z^2-12z-5)\(4\sum_{k=0}^\infty(k+1)f_k^{(4)}z^k+\sum_{k=1}^\infty kf_k^{(4)}z^{k-1}\)
\\&\ -[z^{n+1}](32z^2+8z)\sum_{k=0}^\infty f_k^{(4)}z^k.
\end{align*}

In view of \eqref{kW_k},
\begin{align*}&5(16z^2-1)\sum_{k=1}^\infty\f{kz^{k-1}}{(1+4z)^{k+1}}W_k\l(\f1{1+4z}\r)
\\=&5(16z^2-1)\sum_{m=1}^\infty m(f_m^{(4)}+4s_m)z^{m-1}
\\=&5(16z^2-1)(6+68z+900z^2+\ldots)=-30-340z-4020z^2-\ldots
\end{align*}
Combining this with the final result in the last paragraph, we find that
\begin{align*}&5(16z^2-1)\sum_{k=1}^\infty\f{kz^{k-1}}{(1+4z)^{k+1}}W_k\l(\f1{1+4z}\r)
\\=&(4z+1)(8z-5)\(4\sum_{k=0}^\infty(k+1)f_k^{(4)}z^k+\sum_{k=1}^\infty kf_k^{(4)}z^{k-1}\)
-8z(4z+1)\sum_{k=0}^\infty f_k^{(4)}z^k
\end{align*}
and hence
\begin{align*}&5(4z-1)\sum_{k=1}^\infty\f{kz^{k-1}}{(1+4z)^{k+1}}W_k\l(\f1{1+4z}\r)
\\=&(8z-5)\(4\sum_{k=0}^\infty(k+1)f_k^{(4)}z^k+\sum_{k=1}^\infty kf_k^{(4)}z^{k-1}\)
-8z\sum_{k=0}^\infty f_k^{(4)}z^k.
\end{align*}
This yields the desired \eqref{z-iden}.

The proof of Theorem \ref{auxi} is now complete. \qed

\medskip
\noindent{\it Proof of Theorem \ref{Th1.1}}. In light of Theorem \ref{auxi}, we have
\begin{align*}\sum_{k=0}^\infty\f{45k+8}{40^k}W_k\l(\f 9{10}\r)&=\f{1075}{72}\sum_{k=0}^\infty\f{4k+1}{36^k}f_k^{(4)},
\\\sum_{k=0}^\infty\f{1360k+389}{(-60)^k}W_k\l(\f{16}{15}\r)
&=\f{9225}{32}\sum_{k=0}^\infty\f{4k+1}{(-64)^k}f_k^{(4)},
\\\sum_{k=0}^\infty\f{735k+124}{200^k}W_k\l(\f{49}{50}\r)
&=\f{10125}{784}\sum_{k=0}^\infty\f{60k+11}{196^k}f_k^{(4)},
\\\sum_{k=0}^\infty\f{376380k+69727}{(-320)^k}W_k\l(\f{81}{80}\r)
&=\f{5209600}{243}\sum_{k=0}^\infty\f{17k+3}{(-324)^k}f_k^{(4)},
\\\sum_{k=0}^\infty\f{348840k+47461}{1300^k}W_k\l(\f{324}{325}\r)
&=\f{1314625}{243}\sum_{k=0}^\infty\f{65k+9}{1296^k}f_k^{(4)},
\\\sum_{k=0}^\infty\f{41673840k+4777111}{5780^k}W_k\l(\f{1444}{1445}\r)
&=\f{147758475}{1444}\sum_{k=0}^\infty\f{408k+47}{5776^k}f_k^{(4)}.
\end{align*}
By S. Cooper \cite{Co},
\begin{gather*}\sum_{k=0}^\infty\f{4k+1}{36^k}f_k^{(4)}=\f{6\sqrt{15}}{5\pi},
\ \ \sum_{k=0}^\infty\f{4k+1}{(-64)^k}f_k^{(4)}=\f{32\sqrt{15}}{45\pi},
\\\sum_{k=0}^\infty\f{60k+11}{196^k}f_k^{(4)}=\f{14\sqrt7}{\pi},
\ \ \sum_{k=0}^\infty\f{17k+3}{(-324)^k}f_k^{(4)}=\f{81\sqrt5}{20\pi},
\\\sum_{k=0}^\infty\f{65k+9}{1296^k}f_k^{(4)}=\f{81\sqrt2}{4\pi},
\ \ \sum_{k=0}^\infty\f{408k+47}{5776^k}f_k^{(4)}=\f{76\sqrt{95}}{5\pi}.
\end{gather*}
So we have the desired \eqref{W2}-\eqref{W15}. This concludes the proof. \qed

\section{Congruences related to the identities \eqref{W2}-\eqref{W15}}
 \setcounter{equation}{0}
 \setcounter{conjecture}{0}

In \cite[Section 3]{S13d} the author introduced the polynomials
\begin{equation}\label{Sn(x)}S_n(x)=\sum_{k=0}^n\bi nk^4x^k\ \ (n=0,1,2,\ldots)
\end{equation}
and made conjectures on $\sum_{k=0}^{p-1}S_k(x)$ modulo $p^2$
(with $p$ an odd prime) for each integer $x$ among the numbers
$$1,\,-2,\,\pm4,\,-9,\,12,\,16,\,-20,\,36,\,-64,\,196,\,-324,\,1296,\,5776.$$
See also \cite[Conjectures 49-51]{S19}.

Theorem \ref{Th1.1} and its proof are actually motivated by the following conjecture.
\begin{conjecture}\label{relation} Let $p$ be an odd prime and let $x$
be a $p$-adic integer with $x\not\eq0\pmod p$. Then
\begin{equation}\label{mod p}\sum_{k=0}^{p-1}\f1{(4x)^k}W_k\l(1-\f1x\r)\eq\sum_{k=0}^{p-1}S_k(4x-4)\pmod{p}.
\end{equation}
When
$$x\in\l\{2,\,\pm\f 54,\,\pm4,\,5,\,10,\,-15,\,50,\,-80,\,325,\,1445\r\},$$
we have the further congruence
\begin{equation}\label{mod p^2}\sum_{k=0}^{p-1}\f1{(4x)^k}W_k\l(1-\f1x\r)\eq\sum_{k=0}^{p-1}S_k(4x-4)\pmod{p^2}.
\end{equation}
\end{conjecture}

 The identity \eqref{W2} is motivated by the following conjecture on related congruences.

\begin{conjecture}\label{Conj-W2} {\rm (i)} For any $n\in\Z^+$ we have
$$\f{10^{n-1}}{4n}\sum_{k=0}^{n-1}(45k+8)40^{n-1-k}W_k\l(\f 9{10}\r)\in\Z^+.$$

{\rm (ii)} Let $p\not=2,5$ be a prime. Then
$$\sum_{k=0}^{p-1}\f{45k+8}{40^k}W_k\l(\f 9{10}\r)\eq\f p{16}\l(129\l(\f{-15}p\r)-1\r)\pmod{p^2}.$$
When $(\f{-15}p)=1$, for any $n\in\Z^+$ the number
$$\sum_{k=0}^{pn-1}\f{45k+8}{40^k}W_k\l(\f 9{10}\r)
-p\sum_{k=0}^{n-1}\f{45k+8}{40^k}W_k\l(\f 9{10}\r)$$
divided by $(pn)^2$ is a $p$-adic integer.
\end{conjecture}

The identity \eqref{W3} is motivated by the following conjecture on related congruences.

\begin{conjecture}\label{Conj-W3} {\rm (i)} For any $n\in\Z^+$ we have
$$\f{15^{n-1}}{n}\sum_{k=0}^{n-1}(1360k+389)(-60)^{n-1-k}W_k\l(\f {16}{15}\r)\in\Z^+.$$

{\rm (ii)} Let $p>5$ be a prime. Then
$$\sum_{k=0}^{p-1}\f{1360k+389}{(-60)^k}W_k\l(\f {16}{15}\r)\eq\f p{2}\l(779\l(\f{-15}p\r)-1\r)\pmod{p^2}.$$
When $(\f{-15}p)=1$, for any $n\in\Z^+$ the number
$$\sum_{k=0}^{pn-1}\f{1360k+389}{(-60)^k}W_k\l(\f {16}{15}\r)
-p\sum_{k=0}^{n-1}\f{1360k+389}{(-60)^k}W_k\l(\f {16}{15}\r)$$
divided by $(pn)^2$ is a $p$-adic integer.
\end{conjecture}

The identity \eqref{W6} is motivated by the following conjecture on related congruences.

\begin{conjecture}\label{Conj-W6} {\rm (i)} For any $n\in\Z^+$ we have
$$\f{50^{n-1}}{4n}\sum_{k=0}^{n-1}(735k+124)200^{n-1-k}W_k\l(\f {49}{50}\r)\in\Z^+.$$

{\rm (ii)} Let $p\not=2,5$ be a prime. Then
$$\sum_{k=0}^{p-1}\f{735k+124}{200^k}W_k\l(\f {49}{50}\r)
\eq \f p{32}\l(3969\l(\f{-7}p\r)-1\r)\pmod{p^2}.$$
When $(\f p7)=1$, for any $n\in\Z^+$ the number
$$\sum_{k=0}^{pn-1}\f{735k+124}{200^k}W_k\l(\f {49}{50}\r)
-p\sum_{k=0}^{n-1}\f{735k+124}{200^k}W_k\l(\f {49}{50}\r)$$
divided by $(pn)^2$ is a $p$-adic integer.
\end{conjecture}

The identity \eqref{W8} is motivated by the following conjecture on related congruences.

\begin{conjecture}\label{Conj-W8} {\rm (i)} For any $n\in\Z^+$ we have
$$\f{80^{n-1}}{n}\sum_{k=0}^{n-1}(376380k+69727)(-1)^k320^{n-1-k}W_k\l(\f {81}{80}\r)\in\Z^+.$$

{\rm (ii)} Let $p\not=2,5$ be a prime. Then
$$\sum_{k=0}^{p-1}\f{376380k+69727}{(-320)^k}W_k\l(\f {81}{80}\r)
\eq \f p{3}\l(209198\l(\f{-5}p\r)-17\r)\pmod{p^2}.$$
When $(\f {-5}p)=1$, for any $n\in\Z^+$ the number
$$\sum_{k=0}^{pn-1}\f{376380k+69727}{(-320)^k}W_k\l(\f {81}{80}\r)
-p\sum_{k=0}^{n-1}\f{376380k+69727}{(-320)^k}W_k\l(\f {81}{80}\r)$$
divided by $(pn)^2$ is a $p$-adic integer.
\end{conjecture}

The identity \eqref{W12} is motivated by the following conjecture on related congruences.

\begin{conjecture}\label{Conj-W12} {\rm (i)} For any $n\in\Z^+$ we have
$$\f{325^{n-1}}{n}\sum_{k=0}^{n-1}(348840k+47461)1300^{n-1-k}W_k\l(\f {324}{325}\r)\in\Z^+,$$
and this number is odd if and only if $n\in\{2^a:\ a\in\N\}$.

{\rm (ii)} Let $p\not=2,5,13$ be a prime. Then
$$\sum_{k=0}^{p-1}\f{348840k+47461}{1300^k}W_k\l(\f {324}{325}\r)
\eq \f p{3}\l(142384\l(\f{-2}p\r)-1\r)\pmod{p^2}.$$
When $p\eq1,3\pmod 8$, for any $n\in\Z^+$ the number
$$\sum_{k=0}^{pn-1}\f{348840k+47461}{1300^k}W_k\l(\f {324}{325}\r)
-p\sum_{k=0}^{n-1}\f{348840k+47461}{1300^k}W_k\l(\f {324}{325}\r)$$
divided by $(pn)^2$ is a $p$-adic integer.
\end{conjecture}

The identity \eqref{W15} is motivated by the following conjecture on related congruences.

\begin{conjecture}\label{Conj-W15} {\rm (i)} For any $n\in\Z^+$ we have
$$\f{1445^{n-1}}{n}\sum_{k=0}^{n-1}(41673840k+4777111)5780^{n-1-k}W_k\l(\f {1444}{1445}\r)\in\Z^+,$$
and this number is odd if and only if $n\in\{2^a:\ a\in\N\}$.

{\rm (ii)} Let $p\not=2,5,17$ be a prime. Then
$$\sum_{k=0}^{p-1}\f{41673840k+4777111}{5780^k}W_k\l(\f {1444}{1445}\r)
\eq p\l(4777113\l(\f{-95}p\r)-2\r)\pmod{p^2}.$$
When $(\f{-95}p)=1$, for any $n\in\Z^+$ the number
$$\sum_{k=0}^{pn-1}\f{5928k+253}{5780^k}W_k\l(\f {1156}{5}\r)
-p\sum_{k=0}^{n-1}\f{5928k+253}{5780^k}W_k\l(\f {1156}{5}\r)$$
divided by $(pn)^2$ is a $p$-adic integer.
\end{conjecture}

The conjectural identities \eqref{W1}-\eqref{W14} are motivated by related congruences
stated in \cite[Conjectures 10.34-10.42]{S-book}.

\section{A new type series for $1/\pi$ involving generalized central trinomial coefficients}
\setcounter{equation}{0}
\setcounter{conjecture}{0}

For $b,c\in\Z$ and $n\in\N$ the generalized trinomial coefficient $T_n(b,c)$
denotes the coefficient of $x^n$ in the expansion of $(x^2+bx+c)^n$.

In 2011, the author \cite{S-11,S14c} posed over 60 conjectural series for $1/\pi$ of
the following seven types with $a,b,c,d,m$ integers and $mbcd(b^2-4c)$ nonzero.

\ \ {\tt Type I}. $\sum_{k=0}^\infty\f{a+dk}{m^k}\bi{2k}k^2T_k(b,c)$.

\ \ {\tt Type II}.
$\sum_{k=0}^\infty\f{a+dk}{m^k}\bi{2k}k\bi{3k}kT_k(b,c)$.

\ \ {\tt Type III}.
$\sum_{k=0}^\infty\f{a+dk}{m^k}\bi{4k}{2k}\bi{2k}kT_k(b,c)$.

\ \ {\tt Type IV}.
$\sum_{k=0}^\infty\f{a+dk}{m^k}\bi{2k}{k}^2T_{2k}(b,c)$.

\ \ {\tt Type V}.
$\sum_{k=0}^\infty\f{a+dk}{m^k}\bi{2k}{k}\bi{3k}kT_{3k}(b,c)$.

\ \ {\tt Type VI}.
$\sum_{k=0}^\infty\f{a+dk}{m^k}T_{k}(b,c)^3.$

\ \ {\tt Type VII}.
$\sum_{k=0}^\infty\f{a+dk}{m^k}\bi{2k}kT_{k}(b,c)^2.$
\medskip

Though some of these new families of conjectural series for $1/\pi$
have been proved (see, e.g., \cite{ChWZ}), the three conjectual series for $1/\pi$ of type VI and two of type VII remain open.

In a recent published paper \cite{S20} the author
 proposed four conjectural series for $1/\pi$ of a new type:
\smallskip

\ \ {\tt Type VIII}.
$\sum_{k=0}^\infty\f{a+dk}{m^k}T_k(b,c)T_{k}(b_*,c_*)^2,$
\smallskip
\newline
where $a,b,b_*,c,c_*,d,m$ are integers with $mbb_*cc_*d(b^2-4c)(b_*^2-4c_*)(b^2c_*-b_*^2c)\not=0$.

Here we introduce series for $1/\pi$ involving generalized central trinomial coefficients
of the following novel type:
\smallskip

\ \ {\tt Type IX}.
$\sum_{k=0}^\infty\f{a+dk}{m^k}\bi{2k}kT_k(b,c)T_{k}(b_*,c_*),$
\smallskip
\newline
where $a,b,b_*,c,c_*,d,m$ are integers with $mbb_*cc_*d(b^2-4c)(b_*^2-4c_*)(b^2c_*-b_*^2c)\not=0$.
\medskip

\begin{conjecture}\label{Conj-IX} We have the following identities:
\[\sum_{k=0}^\infty\f{4290k+367}{3136^k}\bi{2k}kT_k(14,1)T_k(17,16)=\f{5390}{\pi}\tag{IX1}\]
and
\[\sum_{k=0}^\infty\f{540k+137}{3136^k}\bi{2k}kT_k(2,81)T_k(14,81)=\f{98}{3\pi}(10+7\sqrt5).\tag{IX2}\]
\end{conjecture}

The conjectural identity (IX1) is motivated by the following conjecture on congruences.

\begin{conjecture}\label{Conj-IX1} {\rm (i)} For any integer $n>1$, we have
\begin{equation*}n\bi{2n}n\ \bigg|\ \sum_{k=0}^{n-1}(4290k+367)3136^{n-1-k}\bi{2k}kT_k(14,1)T_k(17,16).
\end{equation*}

{\rm (ii)} Let $p$ be an odd prime with $p\not=7$. Then
\begin{align*}&\sum_{k=0}^{p-1}\f{4290k+367}{3136^k}\bi{2k}kT_k(14,1)T_k(17,16)
\\\eq&\f p2\l(1430\l(\f{-1}p\r)+39\l(\f 3p\r)-735\r)\pmod{p^2}.
\end{align*}
Moreover, when $p\eq1\pmod{12}$, for any $n\in\Z^+$ the number
\begin{align*}&\sum_{k=0}^{pn-1}\f{4290k+367}{3136^k}\bi{2k}kT_k(14,1)T_k(17,16)
\\&-p\sum_{k=0}^{n-1}\f{4290k+367}{3136^k}\bi{2k}kT_k(14,1)T_k(17,16)
\end{align*}
divided by $(pn)^2\bi{2n}n$ is a $p$-adic integer.

{\rm (iii)} For any prime $p>7$, we have
\begin{align*}&\l(\f{-1}p\r)\sum_{k=0}^{p-1}\f{\bi{2k}k}{3136^k}T_k(14,1)T_k(17,16)
\\\eq&\begin{cases}4x^2-2p\pmod{p^2}&\t{if}\ p\eq1,4\pmod{15}\ \&\ p=x^2+15y^2\,(x,y\in\Z),
\\2p-12x^2\pmod{p^2}&\t{if}\ p\eq2,8\pmod{15}\ \&\ p=3x^2+5y^2\,(x,y\in\Z),
\\0\pmod{p^2}&\t{if}\ (\f{-15}p)=-1.
\end{cases}
\end{align*}
\end{conjecture}
\begin{remark} Note that the imaginary quadratic field $\Q(\sqrt{-15})$ has class number two.
\end{remark}

The conjectural identity (IX2) is motivated by the following conjecture on congruences.

\begin{conjecture}\label{Conj-IX2} {\rm (i)} For any integer $n>1$, we have
\begin{equation*}2n\bi{2n}n\ \bigg|\ \sum_{k=0}^{n-1}(540k+137)3136^{n-1-k}\bi{2k}kT_k(2,81)T_k(14,81).
\end{equation*}

{\rm (ii)} Let $p$ be an odd prime with $p\not=7$. Then
\begin{align*}&\sum_{k=0}^{p-1}\f{540k+137}{3136^k}\bi{2k}kT_k(2,81)T_k(14,81)
\\\eq&\f p3\l(270\l(\f{-1}p\r)-104\l(\f {-2}p\r)+245\l(\f{-5}p\r)\r)\pmod{p^2}.
\end{align*}
Moreover, when $p\eq\pm1,\pm9\pmod{40}$, for any $n\in\Z^+$ the number
\begin{align*}&\sum_{k=0}^{pn-1}\f{540k+137}{3136^k}\bi{2k}kT_k(2,81)T_k(14,81)
\\&-p\l(\f{-1}p\r)\sum_{k=0}^{n-1}\f{540k+137}{3136^k}\bi{2k}kT_k(2,81)T_k(14,81)
\end{align*}
divided by $(pn)^2\bi{2n}n$ is a $p$-adic integer.

{\rm (iii)} For any prime $p>7$, we have
\begin{align*}&\l(\f{-1}p\r)\sum_{k=0}^{p-1}\f{\bi{2k}k}{3136^k}T_k(2,81)T_k(14,81)
\\\eq&\begin{cases}4x^2-2p\pmod{p^2}&\t{if}\ (\f 2p)=(\f p3)=(\f p5)=1\ \&\ p=x^2+30y^2,
\\8x^2-2p\pmod{p^2}&\t{if}\ (\f 2p)=1,\ (\f p3)=(\f p5)=-1\ \&\ p=2x^2+15y^2,
\\20x^2-2p\pmod{p^2}&\t{if}\ (\f p5)=1,\ (\f 2p)=(\f p3)=-1\ \&\ p=5x^2+6y^2,
\\2p-12x^2\pmod{p^2}&\t{if}\ (\f p3)=1,\ (\f 2p)=(\f p5)=-1\ \&\ p=3x^2+10y^2,
\\0\pmod{p^2}&\t{if}\ (\f{-30}p)=-1,
\end{cases}
\end{align*}
where $x$ and $y$ are integers.
\end{conjecture}
\begin{remark}\label{Rem5.6} Note that the imaginary quadratic field $\Q(\sqrt{-30})$
has class number four.
\end{remark}

\section{Other new conjectural series for $1/\pi$}
\setcounter{equation}{0}
\setcounter{conjecture}{0}

As mentioned in \cite[Remark 4.4]{S14d}, an identity of MacMahon implies that the polynomial
$$F_n(x)=\sum_{k=0}^n\bi nk\bi{n+2k}{2k}\bi{2k}kx^{n-k}$$
at $x=-4$ coincides with the Franel number $f_n=\sum_{k=0}^n\bi nk^3$.
Conjecture 4.4 of Sun \cite{S14d} lists ten conjectural series for $1/\pi$
involving $F_n(x)$ with $x\not=-4$; eight of them were later confirmed in \cite{CWZ}, but the following two remain open:
\begin{align}\label{5.1}\sum_{k=0}^\infty\f{357k+103}{2160^k}\bi{2k}kF_k(-324)=&\f{90}{\pi},
\\\label{5.2} \sum_{k=0}^\infty\f k{3645^k}\bi{2k}kF_k(486)=&\f{10}{3\pi}.
\end{align}

Here we pose the following new conjecture.

\begin{conjecture}\label{Conj5.1} We have the following identities:
\begin{align}\label{5.3}\sum_{k=0}^\infty\f{6k+1}{(-1728)^k}\bi{2k}kF_k(-324)&=\f{24}{25\pi}\sqrt{375+120\sqrt{10}},
\\\label{5.4}\sum_{k=0}^\infty\f{4k+1}{(-160)^k}\bi{2k}kF_k(-20)&=\f{\sqrt{30}}{5\pi}\cdot\f{5+\root3\of{145+30\sqrt6}}
{\root6\of{145+30\sqrt6}},
\\\label{5.5}\sum_{k=0}^\infty\f{1290k+289}{27648^k}\bi{2k}kF_k(-2160)&=\f{96\sqrt{15}}{\pi},
\\\label{5.6}\sum_{k=0}^\infty\f{804k+49}{276480^k}\bi{2k}kF_k(12096)&=\f{120\sqrt{15}}{\pi},
\\\label{5.7}\sum_{k=0}^\infty(24k+5)\l(\f 2{135}\r)^kF_k\l(-\f{27}8\r)&=\f{3}{2\pi}(5\sqrt6+4\sqrt{15}).
\end{align}
\end{conjecture}
\begin{remark} The author actually found $(5.3)$-$(5.7)$ in 2020.
As all of them converge quickly, one can easily check them via {\tt Mathematica} or {\tt Maple}.
\end{remark}

The identity \eqref{5.3} is motivated by \cite[Conjecture 4.6]{S14d}. The reader might wonder how we found the right-hand side
of the identity \eqref{5.3}. We thought that the left-hand side of \eqref{5.3} times $\pi$
is an algebraic number and found the form of this algebraic number via
calculating its first 100 digits and using the Maple command {\tt identify}.

The identities \eqref{5.4} and \eqref{5.5} are motivated by related congruences stated in
\cite[Conjectures 10.47--10.48]{S-book}.

The identity \eqref{5.6} is motivated by the following conjecture on related congruences.

\begin{conjecture}\label{Conj5.5} {\rm (i)} Let $n>1$ be an integer. Then
$$\f1{n\bi{2n}n}\sum_{k=0}^{n-1}(804k+49)276480^{n-1-k}\bi{2k}kF_k(12096)\in\Z^+,$$
and this number is odd if and only if $n\in\{2^a+1:\ a\in\N\}$.

{\rm (ii)} Let $p>5$ be a prime. Then
\begin{align*}\sum_{k=0}^{p-1}\f{804k+49}{276480^k}\bi{2k}kF_k(12096)
\eq p\l(95\l(\f {-15}p\r)-46\l(\f{30}p\r)\r)\pmod{p^2}.
\end{align*}
Moreover, if $p\eq1,3\pmod8$ then for any $n\in\Z^+$ the number
$$\sum_{k=0}^{pn-1}\f{804k+49}{276480^k}\bi{2k}kF_k(12096)-p\l(\f {-15}p\r)
\sum_{k=0}^{n-1}\f{804k+49}{276480^k}\bi{2k}kF_k(12096)$$
divided by $(pn)^2\bi{2n}n$ is a $p$-adic integer.

{\rm (iii)} Let $p>5$ be a prime. Then
\begin{align*}&\sum_{k=0}^{p-1}\f{\bi{2k}k}{276480^k}F_k(12096)
\\\eq&\begin{cases}4x^2-2p\pmod{p^2}&\t{if}\ (\f{-2}p)=(\f p3)=(\f p5)=(\f p7)=1\ \&\ p=x^2+210y^2,
\\8x^2-2p\pmod{p^2}&\t{if}\ (\f{-2}p)=(\f p7)=1,\ (\f p3)=(\f p5)=-1\ \&\ p=2x^2+105y^2,
\\2p-12x^2\pmod{p^2}&\t{if}\ (\f{-2}p)=(\f p3)=1,\ (\f p5)=(\f p7)=-1\ \&\ p=3x^2+70y^2,
\\20x^2-2p\pmod{p^2}&\t{if}\ (\f{-2}p)=(\f p3)=(\f p5)=(\f p7)=-1\ \&\ p=5x^2+42y^2,
\\2p-24x^2\pmod{p^2}&\t{if}\ (\f{-2}p)=(\f p5)=1,\ (\f p3)=(\f p7)=-1\ \&\ p=6x^2+35y^2,
\\28x^2-2p\pmod{p^2}&\t{if}\ (\f{-2}p)=(\f p5)=-1,\ (\f p3)=(\f p7)=1\ \&\ p=7x^2+30y^2,
\\40x^2-2p\pmod{p^2}&\t{if}\ (\f{-2}p)=(\f p7)=-1,\ (\f p3)=(\f p5)=1\ \&\ p=10x^2+21y^2,
\\56x^2-2p\pmod{p^2}&\t{if}\ (\f{-2}p)=(\f p3)=-1,\ (\f p5)=(\f p7)=1\ \&\ p=14x^2+15y^2,
\\0\pmod{p^2}&\t{if}\ (\f{-210}p)=-1,
\end{cases}
\end{align*}
where $x$ and $y$ are integers.
\end{conjecture}
\begin{remark}\label{Rem5.5} Note that the imaginary quadratic field $\Q(\sqrt{-210})$
has class number eight.
\end{remark}

The identity \eqref{5.7} is motivated by the following conjecture on related congruences.

\begin{conjecture}\label{Conj2.6} {\rm (i)} Let $n$ be any positive integer. Then
$$\f{4^{n-1}}{n\bi{2n-1}{n-1}}\sum_{k=0}^{n-1}(24k+5)135^{n-1-k}2^k\bi{2k}kF_k\l(-\f{27}8\r)\in\Z^+,$$
and this number is congruent to $5$ modulo $8$.

{\rm (ii)} Let $p>5$ be a prime. Then
\begin{align*}\sum_{k=0}^{p-1}\f{(24k+5)2^k}{135^k}\bi{2k}kF_k\l(-\f{27}8\r)
\eq p\l(4\l(\f {-6}p\r)+\l(\f{-15}p\r)\r)\pmod{p^2}.
\end{align*}
Moreover, if $(\f{10}p)=1$ then for any $n\in\Z^+$ the number
$$\sum_{k=0}^{pn-1}\f{(24k+5)2^k}{135^k}\bi{2k}kF_k\l(-\f{27}8\r)-p\l(\f {-6}p\r)
\sum_{k=0}^{n-1}\f{(24k+5)2^k}{135^k}\bi{2k}kF_k\l(-\f{27}8\r)$$
divided by $(pn)^2\bi{2n}n$ is a $p$-adic integer.

{\rm (iii)} Let $p>5$ be a prime. Then
\begin{align*}&\sum_{k=0}^{p-1}\f{2^k\bi{2k}k}{135^k}F_k\l(-\f{27}8\r)
\\\eq&\begin{cases}4x^2-2p\pmod{p^2}&\t{if}\ (\f{2}p)=(\f p3)=(\f p5)=1\ \&\ p=x^2+30y^2,
\\8x^2-2p\pmod{p^2}&\t{if}\ (\f{2}p)=1,\ (\f p3)=(\f p5)=-1\ \&\ p=2x^2+15y^2,
\\2p-12x^2\pmod{p^2}&\t{if}\ (\f p3)=1,\ (\f 2p)=(\f p5)=-1\ \&\ p=3x^2+10y^2,
\\20x^2-2p\pmod{p^2}&\t{if}\ (\f{p}5)=1,\ (\f 2p)=(\f p3)=-1\ \&\ p=5x^2+6y^2,
\\0\pmod{p^2}&\t{if}\ (\f{-30}p)=-1,
\end{cases}
\end{align*}
where $x$ and $y$ are integers.
\end{conjecture}

In 2012 the author (cf. \cite[(8)]{S14c}) conjectured that
\begin{equation*}\label{1.1}\sum_{n=0}^\infty\f{28n+5}{576^n}\bi{2n}n\sum_{k=0}^n\f{5^k\bi{2k}k^2\bi{2(n-k)}{n-k}^2}{\bi nk}=\f 9{\pi}(2+\sqrt2),
\end{equation*}
which remains open up to now. Here we pose a similar conjecture.

\begin{conjecture}\label{Conj6.1} We have the following identity:
\begin{equation}\label{1.2}
\sum_{n=0}^\infty\f{182n+31}{576^n}\bi{2n}n\sum_{k=0}^n\f{\bi{2k}k^2\bi{2(n-k)}{n-k}^2}{\bi nk}\l(-\f{25}{16}\r)^k=\f {189}{2\pi}.
\end{equation}
\end{conjecture}

This is motivated by the author's following conjecture on related congruences.
\begin{conjecture}\label{Conj6.2} Let $p>3$ be a prime. Then
\begin{align*}&\sum_{n=0}^{p-1}\f{182n+31}{576^n}\bi{2n}n\sum_{k=0}^n\f{\bi{2k}k^2\bi{2n-2k}{n-k}^2}{\bi nk}\l(-\f{25}{16}\r)^k
\\&\quad\eq\f p2\l(63\l(\f{-1}p\r)-1\r)\pmod{p^2}.
\end{align*}
Also,
\begin{align*}&\sum_{n=0}^{p-1}\f{\bi{2n}n}{576^n}\sum_{k=0}^n\f{\bi{2k}k^2\bi{2n-2k}{n-k}^2}{\bi nk}\l(-\f{25}{16}\r)^k
\\\eq&\begin{cases}4x^2-2p\pmod{p^2}&\t{if}\ (\f p7)=1\ \&\ p=x^2+7y^2\, (x,y\in\Z),
\\0\pmod{p^2}&\t{if}\ (\f p7)=-1,\ \t{i.e.},\ p\eq3,5,6\pmod 7.
\end{cases}\end{align*}
\end{conjecture}

\subsection*{Acknowledgement}
 The work was supported by the National Natural Science Foundation of China (Grant No. 11971222).

\end{document}